\documentclass[11 pt,a4paper,twoside,reqno]{amsart} 
\usepackage{amsfonts,amssymb,amscd,amsmath,enumerate,verbatim,calc} 
 
\begin{document}

\begin{abstract}
 For natural numbers $n$ and $d$, set $p(n)$ to be the number of integer partitions of $n$. Set $a(n)$ to be the number of integer partitions of $n$, all of whose parts fail to be divisible by d. Set $b(n)$ to be the number of integer partitions of n, where no part is repeated d or more times. ‘Glaisher’s theorem
gives a bijective proof of the fact that $a(n) = b(n)$. We present a new family of bijections to show that $p(n)-a(n) = p(n)-b(n)$, which is equivalent to Glaisher's theorem, and is a generalization of Glaisher's original argument. This provides a rich family of elementary bijections between these two sets of partitions.

\end{abstract}

\title{Glaisher's partition problem}
\author{Aritro Pathak,\\  Department of Mathematics \\ Brandeis University, Waltham, MA, 02453, USA}

\maketitle

Glaisher's theorem states that the number of partitions of an integer $N$ into parts not divisible by $d$ is equal to the number of partitions where no part is repeated $d$ or more times.

 This theorem can also be proven in an elementary way using generating functions. The original combinatorial proof is well known, as can be found in \cite{lehmer}. It also appears as Exercise 3.2.3 in Igor Pak's survey on Partition Bijections \cite{pak}. We present an outline of Glaisher's original proof in our first proof below. In the second proof, we give a new family of bijections generalizing Glaisher's argument.
 
 \textbf{Theorem:} Given integers $n\geq d\geq 1$, the number of integer partitions of $n$ all of whose parts fail to be divisible by $d$, is equal to the number of integer partitions of $n$ where no part is repeated $d$ or more times.

\textbf{Proof 1:}  For a partition where each part appears less than $d$ times, split up the parts divisible by $d^{t}$ into $d^{t}$ parts (where $t\geq 0$ is the highest power of $d$ that divides that part). This gives a partition with no part divisible by $d$.

For the inverse, for any partition where no part is divisible by $d$, write the number of times $s^{(j)}$ appears in the base $d$ representation: $\sum_{i} a^{(j)}_{i}d^{i}$. Then consider the partition where $d^{i}s^{(j)}$ each appear $a^{(j)}_{i}<d$ number of times.

It can be checked that this gives a bijection. $ \square$
 
By considering the complement of the two sets under consideration, we construct the new family of bijections in the next proof. This obviously implies Glaisher's theorem.

\textbf{Proof 2:}  We show that the number of elements in the set $A_{N}$ that are partitions of $N$ into parts where at least one part is divisible by $d$, is equal to the number of elements in the set $B_{N}$ that are  partitions where at least one part is repeated $d$ or more times. 

Consider a partition $P_{A}$ in $A_{N}$. For each $x\in \mathbb{Z}_{+}$ that is not divisible by $d$, construct a square matrix $M^{(x)}_{ij}, \ (i,j\in \mathbb{Z}_{+})$, where the square at the intersection of the row $i$ and column $j$ contains the coefficient of $d^{i}$ in the base $d$ expansion of the integer that equals the number of times that $x\cdot d^{j}$ appears in $P_{A}$. 

By the definition of $A_{N}$, for this partition, there exists at least one $x_{P_{A}}\in \mathbb{Z}_{+}$ not divisible by $d$, so that the matrix $M^{(x_{P_{A}})}$ has a non zero element in the columns indexed by $j\geq 2$ (i.e it's not just the first column that is non empty).

Similarly, for any given partition $P_{B}$ in $B_{N}$, we can construct for each $x\in \mathbb{Z}_{+}$ not divisible by $d$, a square matrix $\Tilde{M}^{(x)}_{ij}$ as above indexed by $\mathbb{Z}_{+}\times \mathbb{Z}_{+}$ so that there exists some $x_{P_{B}}$ not divisible by $d$, whose matrix contains a non zero element  in the rows corresponding to $i\geq 2$ (i.e it's not just the first row that has non zero elements). 

Now we construct the family of bijections between $A_{N}$ and $B_{N}$. 

For the partition $P_{A}$ in $A_{N}$, looking at the matrix corresponding to $x_{P_{A}}$, we can look at the "southwest-northeast diagonals" indexed by $k\in \mathbb{Z}_{+}$: $D^{(k)}_{x_{P_{A}}}=\{(i,j)\in M_{ij}^{(x_{P_{A}})}|i+j=k, i,j\in \mathbb{Z}_{+}\}$.  

In any diagonal indexed by some integer $k$ , we can take a permutation with the restriction that the $(1,k-1)$ square is taken in the interior of this diagonal, which is the set $\text{int}(D^{(k)}_{x_{P_{A}}})=\{(i,j)| (i,j)\in  D^{(k)}_{x_{P_{A}}}, i\neq k-1 , \ j \neq k-1 \}$, and the $(k-1,1)$ square is taken to the $(1, k-1)$ square. 

Since we are permuting within each fixed diagonal, the contribution to the sum of $N$ is invariant by this bijection.

It is clear that this bijectively gives us an element $\Tilde{P}_{A}\in B_{N}$, and on applying the corresponding inverse permutation to all the diagonals corresponding to the matrix in $\Tilde{P}_{A}\in B_{N}$, we recover $P_{A}\in A_{N}$.

Note that it is important in general, that in the previous permutation in going from $P_{A}\in A_{N}$ to $\Tilde{P_{A}}\in B_{N}$, we don't take squares from $\text{int}(D^{(k)}_{x_{P_{A}}})$ to go the upper row consisting of the squares $\{(1,k), k\geq 2\}$ by the permutation. This is to negate the cases where the only non zero elements in $M^{(x_{P_{A}})}$ are precisely those in the interior of the matrix, $\cup_{k=2}^{\infty} \text{int}(D^{(k)}_{x_{P_{A}}}) $, that get permuted to go into the upper row, and hence we get a matrix where all the rows corresponding to $i\geq 2$ are empty, so we do not get an element of $B_{N}$. Sending the $(k-1,1)$ square in $M^{x_{P_{A}}}_{ij}$ to go to $(1,k-1)$ in ${\Tilde{M}^{x_{\Tilde{P}_{A}}}_{ij}}$ ensures we have a well defined bijection from $A_{N}$ to $B_{N}$.

This establishes the bijection between partitions in $A_{N}$ and partitions in $B_{N}$. $\square$
\\
\textbf{Remark:} The simplest possible permutation of any diagonal $D_{k}$ is to simply interchange $(1,k-1)$ and $(k-1,1)$, while keeping all the elements of $\text{int}(D^{(k)}_{x_{P_{A}}})$ fixed or choosing a random permutation restricted to $\text{int}(D^{(k)}_{x_{P_{A}}})$. However, the proof above gives a bigger class of permutations.

In Glaisher's original argument for the corresponding complement sets, in one partition set only the first row of our matrix has non zero elements, while in the other partition set, only the first column has non zero elements. Glaisher's bijection involved simply swapping the $(1,j)$ th square with the $(j,1)$ square. So the simplest permutations of the previous paragraph, where the ends of the diagonal and the interior of the diagonal are permuted disjointedly, are `Glaisher-like', while our permutations where one end of the diagonal is taken inside the interior is a variant of Glaisher's original argument.

\textbf{Acknowledgement}: The author is thankful to Konstantin Matveev for pointing out this problem in a Combinatorics class at Brandeis University. The author is also thankful to Dmitry Kleinbock and Ira Gessel for their suggestions upon looking at the manuscript.

\end{document}